\documentclass[12pt]{amsart}
\usepackage{amssymb}
\usepackage{amsmath}
\usepackage{latexsym}
\usepackage{amsthm}
\usepackage{rotate}
\usepackage{graphicx}
\usepackage[normalem]{ulem}
\usepackage{cancel}
\usepackage[utf8]{inputenc}

\usepackage{tikz}
    \usetikzlibrary{arrows,topaths}

 \makeatletter \textwidth 6.2in \textheight
8.6in \oddsidemargin 0.0in \evensidemargin 0.0in \topmargin
0.0in
\def\thmhead@plain#1#2#3{%
 \thmname{#1}\thmnumber{\@ifnotempty{#1}{
 }#2}%
 \thmnote{ \the\thm@notefont(#3)}}
\let\thmhead\thmhead@plain
\def\swappedhead#1#2#3{%
 \thmnumber{#2}\thmname{\@ifnotempty{#2}{. }#1}%
 \thmnote{ \the\thm@notefont(#3)}}
\makeatother

\theoremstyle{definition} 
\newtheorem{definition}{Definition}[section]
\newtheorem{remark}[definition]{Remark}

\theoremstyle{plain}      
\newtheorem{proposition}[definition]{Proposition}
\newtheorem{theorem}[definition]{Theorem}
\newtheorem{corollary}[definition]{Corollary}
\newtheorem{lemma}[definition]{Lemma}

\begin{document}
\email{aa145@aub.edu.lb}

\title{Non-binary counterparts of the Prefer-same and Prefer-opposite De Bruijn Sequences}
\author[A. Alhakim]{Abbas Alhakim\\ Department of Mathematics\\ American University of Beirut\\ Beirut, Lebanon}
\begin{abstract}
The well-known prefer-one, prefer-opposite, and prefer-same binary
de~Bruijn sequences are all constructed using simple preference
rules. We apply the technique of preference functions of span one to
define $q$-ary sequences that generalize the prefer-opposite and
prefer-same sequences and we present some of their basic properties
that are shared with their binary versions. In particular, we show
that the prefer-higher sequence (the nonbinary counter-part of the
prefer-one sequence) is obtained from a homomorphic image of the
proposed prefer-opposite, when repetitions are cleaned up. This
mirrors a known relationship between the binary versions. We also
perform calculations that demonstrate that the discrepancy profile
of the proposed sequences is similar to that of the binary case.
\end{abstract}
\footnotetext[1]{This research was partially supported by the
University Research Board (URB) of the American University of
Beirut. Project Number 104107} \maketitle

\section{Introduction}
Given positive integers $q>1$ and $n\geq1$, a $q$-ary de Bruijn
sequence of order $n$ is a sequence of symbols such that every
pattern of size $n$ appears exactly once as a contiguous
subsequence. %
A small class of particularly simple greedy algorithms that generate
binary de~Bruijn sequences of any given order are the prefer-one
algorithm, see \cite{Ford1957, Martin1934}, the prefer-same
algorithm~\cite{Eldert, Fred1982} and the prefer-opposite algorithm
~\cite{Alhakim10}.

The oldest and most well-known of these is the prefer-one. This
algorithm proceeds by appending bits to an initial string of $n$
zeros in such a way that every time a bit is to be added, $1$ is
attempted first and appended to the sequence. If the last $n$ bits
form a word that was not encountered earlier then this appended $1$
is kept. Otherwise it is replaced by a $0$ and the same check is
made. The process is repeated in the same fashion and it is stopped
as soon as neither 1 nor 0 puts a new $n$-tuple. The prefer-opposite
algorithm works in a similar way but favors the opposite of the most
recently placed bit. Likewise, the prefer same algorithm first
proposes to duplicate the most recently placed bit. In its original
version, it must be started with a string of $n$ ones followed by a
string of $n$ zeros and the algorithm must also keep track of the
number of occurrences of runs of consecutive ones and zeros of all
possible lengths in order to know whether or not a proposed bit is
to be accepted, even when it puts a new $n$-tuple. \cite{Alhakim19}
shows that the extra counting of runs is not necessary, provided the
initial word of the sequence is the alternating string $1010\cdots$,
thus establishing that the prefer-same algorithm is indeed of the
same type as its other two sequences.

Another attractive feature of the three sequences is that they are
all extremes. It is well known that the Prefer-one sequence is the
lexicographically least de Bruijn sequence for a given order $n$. On
the other hand, the prefer-opposite and prefer-same sequences
respectively have the smallest and lexicographically largest run
length encoding, see~\cite{Alhakim19} for definition and proof.


When $q$ is higher than $2$, a generalization of the prefer-one
sequence, known as the prefer-higher or non-binary Ford sequence, is
classical, see Fredricksen~\cite{Fred1982}. For the other two
sequences, the concept of `opposite' is needed but not well defined.
We will employ preference functions of span $1$ to obtain natural
generalization of the prefer-opposite and prefer-same sequences to
the case of alphabet size $q$. Span-1 preference functions
correspond to those de Bruijn sequences that are representable by
$q\times q$ matrices. Indeed, the binary de Bruijn sequences
mentioned above are representable using $2\times2$ matrices together
with an appropriate starting pattern.


The rest of the paper is organized as follows. In
Section~\ref{S:prelims} we give basic definitions and background
about preference functions, with preliminary lemmas that will be
essential for the rest of the paper. In Section 3 we propose two
classes of nonbinary sequences which we respectively call nonbinary
prefer-opposite and prefer-same sequences and we establish
properties that are similar to those of the binary counterparts. In
Section 4 we show that mapping the proposed prefer-opposite sequence
via an appropriate homomorphism we recover the prefer-higher
sequence, replicating the relationship between the binary
prefer-opposite and the prefer-one sequences, via the Lempel
Homomorphism.

\section{Preliminaries}\label{S:prelims}
In this section we give basic definitions and provide important
lemmas that are essential for establishing the results of this
paper. For $q\geq2$, $\mathcal{A}$ is an alphabet with $q$ symbols
that we will denote as $\{0,1,\ldots,q-1\}$. We will interchangeably
refer to these symbols as digits. For an integer $n\geq1$, a
sequence of $n$ digits will be referred to as an $n$-\emph{word},
and denoted as $a_1\cdots a_n$ and often as $(a_1,\ldots,a_n)$ for
notational clarity. $a^n$ denotes the word obtained by repeating the
symbol $a$ $n$ times. A \emph{de Bruijn digraph} of order $n$ admits
all possible $n$-words as \emph{vertices}, and for two vertices
${\bf{a}}=(a_1,\ldots,a_n)$ and ${\bf{b}}=(b_1,\ldots,b_n)$, there
is an \emph{edge} from $\bf{a}$ to  $\bf{b}$   if and only if
$a_i=b_{i-1}$ for $i=2,\ldots,n$. A string $a_1\ldots a_L$ of length
$L\geq n$ can then be identified with a \emph{directed path} (or
simply a path) in the de~Bruijn digraph with vertices
$(a_1,\ldots,a_n), (a_2,\ldots,a_{n+1}),\cdots,
(a_{L-n+1},\ldots,a_L)$. A \emph{cycle} is a path whose first and
last vertices coincide. A \emph{simple cycle} is one that does not
cross a vertex, other than the first and last, more than once. The
first definition below follows Golomb \cite{Golomb1967}.
\begin{definition}
For an integer $n\geq1$, a preference function is a function $P:
\mathcal{A}^{n}\rightarrow S$, where $S$ is the set of all
permutations of the elements of $\mathcal{A}$. We write
$P(\mathbf{x})=(P_1,\ldots,P_q)$ for every $n$-word $\bf{x}$.
\end{definition}

\begin{definition}
For an integer $s,\;0\leq s\leq n$, a preference function $P$
defined on $\mathcal{A}^n$ is said to be of span $s$ if $s$ is the
smallest integer such that $P(x_1,\ldots,x_n)$ is fully determined
by $(x_{n-s+1},\ldots,x_n$, for all $n$-words $(x_1,\ldots,x_n)$.
\end{definition}

Examples of preference functions of spans $0,1,2$ for $q=3$ are
displayed as diagrams in Table~\ref{Ta:Pref_Diagrams}, where the
permutation to the right of an arrow is the function value of the
pattern on the left of the arrow. Observe that the upper left
diagram is that of a constant preference function and so it is of
span $0$.


\begin{table}
\small
\centering
\begin{tabular}{c|c|c}
\hline
\begin{tabular}{lllll}
0 & $\rightarrow$ & 2, & 1, & 0\\
1 & $\rightarrow$ & 2, & 1, & 0\\
2 & $\rightarrow$ & 2, & 1, & 0\\
\end{tabular}
&
\begin{tabular}{lllll}
0 & $\rightarrow$ & 1, & 2, & 0\\
1 & $\rightarrow$ & 1, & 0, & 2\\
2 & $\rightarrow$ & 2, & 1, & 0\\
\end{tabular}
&
\begin{tabular}{lllll}
0 & $\rightarrow$ & 1, & 2, & 0\\
1 & $\rightarrow$ & 2, & 1, & 0\\
2 & $\rightarrow$ & 0, & 1, & 2\\
\end{tabular}\\\hline\hline
\begin{tabular}{lllll}
00 & $\rightarrow$ & 0, & 2, & 1\\
01 & $\rightarrow$ & 1, & 2, & 0\\
02 & $\rightarrow$ & 2, & 1, & 0\\
\end{tabular}
&
\begin{tabular}{lllll}
10 & $\rightarrow$ & 1, & 2, & 0\\
11 & $\rightarrow$ & 0, & 1, & 2\\
12 & $\rightarrow$ & 0, & 2, & 1\\
\end{tabular}
&
\begin{tabular}{lllll}
20 & $\rightarrow$ & 1, & 2, & 0\\
21 & $\rightarrow$ & 0, & 2, & 1\\
22 & $\rightarrow$ & 0, & 1, & 2\\
\end{tabular}
\\\hline
\end{tabular}
\caption{Top row: three preference diagrams of span 1. Lower row:
one preference diagram of span 2.}\label{Ta:Pref_Diagrams}
\end{table}


The following recursive construction produces a unique finite
sequence $\{a_i\}$ from an alphabet $\mathcal{A}$, given an
arbitrary initial $n$-word $(I_1,\cdots,I_n)$ with $n>s$ and a
preference function of span $s$. We denote the resulting sequence by
$(P,I)$.

1. For $i=1,\cdots,n$ let $a_i=I_i$.

2. Suppose that $a_1,\cdots,a_k$ for some integer $k\geq n$ have
been defined. Let $a_{k+1}= P_i(a_{k-s+1},\cdots,a_k)$ where $i,\;1\leq i\leq q$ is
the smallest integer such that the $n$-word
$(a_{k-n+2},\cdots,a_{k+1})$ has not appeared in the sequence as a
substring, if such an $i$ exists.

3. If no such $i$ exists, halt the program (the construction is
complete.)

The following lemma is a slight generalization of Lemma~2 of
Chapter~3 in Golomb~\cite{Golomb1967}. The proof is essentially the
same.
\begin{lemma}\label{L:Golomb}
Consider an arbitrary preference function $P$ of span $s\geq0$ and
initial word $I=(I_1,\cdots,I_n); n>s$. Then every $n$-word occurs
at most once in $(P,I)$. Furthermore, $(P,I)$ ends with the pattern
$(I_1,\cdots,I_{n-1})$.
\end{lemma}

It follows that the sequence $(P,I)$ can be identified with a simple
cycle in the de~Bruijn digraph of order $n$, by removing the last
pattern $(I_1,\ldots,I_{n-1})$ and wrapping the rest around a
circle. In the following two definitions, we fix an integer $s\geq1$
and a preference function $P$ of span $s$, which can equivalently be
considered as a function from $\mathcal{A}^s$ to $S$.
\begin{definition}
For an integer $i$ such that $1\leq i\leq q$, the $i^{th}$ column
function induced by $P$ is a function from $\mathcal{A}^s$ to
$\mathcal{A}^s$ defined as
\[
g_i(x_1,\cdots,x_s)=(x_2,\cdots,x_s,P_i(x_1,\cdots,x_s)).
\]
\end{definition}

Naturally, $g_i$ will define at least one cycle of length $k\geq1$.
That is, a sequence of $k$ vertices $v_1,\cdots,v_k$ in
$\mathcal{A}^s$ such that $g_i(v_j)=v_{j+1}$ for $j=1,\cdots,k-1$
and $g_i(v_k)=v_1$.

\begin{definition}~\label{D:word_on_cycle}
For a given $i$, $1\leq i\leq q$ and a cycle $C$ defined by $g_i$,
we say that a word $(x_1,\cdots,x_s)$ is on $C$ if it is simply one
of the vertices on $C$. For an arbitrary integer $n>s$, we say that
$(x_1,\cdots,x_n)$ is on $C$ if it forms a path in $C$. That is, if
$(x_1,\cdots,x_s)$ is a vertex on $C$ and
$x_j=P_i(x_{j-s},\cdots,x_{j-1})$ for $j=s+1,\cdots,n$.
\end{definition}

We note that a path in $C$ in the above definition need not be a
simple path. That is, it may cross itself one or more times,
depending on the value of $n$. As an example, the word $0^n$ of $n$ zeros is a path on the cycle that sends $0^s$ to itself, for $n>s$.

\begin{definition}\label{D:AssocSets}
Let $g_k$, for some $k$, $1\leq k\leq q$, be the function induced by
the $k^{th}$ column of a preference function $P$. To a given cycle
$C$ on $g_k$ we associate the following sets.

$\bullet$ $[C]$ is the set of all vertices on $C$.

$\bullet$ The closure $\overline{[C]}$ consists of every vertex in
$\mathcal{A}^s$ such that $g_k$ defines a path that starts with this
vertex and terminates in $C$. That is,
$(a_1,\cdots,a_s)\in\overline{[C]}$ if and only if there is a
minimal integer $s^{\prime}\geq s$ such that for
$i=1,\cdots,{s^{\prime}-s}$,
$g_k(a_i,\cdots,a_{i+s-1})=(a_{i+1},\cdots,a_{i+s})$ and
$(a_{s^{\prime}-s+1},\cdots,a_{s^{\prime}})\in[C]$.

$\bullet$ $\Sigma$ is the complement of $\overline{[C]}$ in
$\mathcal{A}^s$.
\end{definition}

The following theorem and its corollary are stated and proved in
Alhakim~\cite{Alhakim21}.

\begin{theorem}\label{T:main}
Let $P$ be a preference function of span $s\geq0$ such that the
induced $g_q$ has [at least] a cycle $C$, and let $[C]$,
$\overline{[C]}$ and $\Sigma$ be the sets associated with $C$ as in
Definition~\ref{D:AssocSets}. For an arbitrary integer  $n>s$ let
$I=(I_1,\cdots,I_n)$ be such that $(I_1,\cdots,I_{n-1})$ is on $C$
as in Definition~\ref{D:word_on_cycle}.

Also, assume that there is an integer $q^{\prime}; 1<q^{\prime}\leq
q$ with the property that $\Sigma$ is closed under all functions
$g_k$ with $k\geq q^{\prime}$, in the sense that whenever
$(x_1,\cdots,x_s)\in\Sigma$, then $g_k(x_1,\cdots,x_s)\in\Sigma$
for all $k\geq q^{\prime}$. Then the sequence $(P,I)$ misses the set
of words $M$ defined as
\[\{(x_1,\cdots,x_n): \left((x_1\cdots x_s)\in\Sigma\right)\;\wedge\;\left(\forall i>s,  x_i=P_k(x_{i-s},\cdots,x_{i-1}) \text{ for some }
k\geq q^{\prime}\right)\}
\]
Furthermore, if $q^{\prime}$ is the smallest integer with the given
property and if the function $\tilde{g}_{q^{\prime}-1}$ (defined as
$g_{q^{\prime}-1}$ restricted to $\Sigma$) has no cycles, then no
other word is missing.
\end{theorem}

\begin{corollary}\label{C:complete}
Let $P$ be a  preference function of span $s\geq0$ such that  the
induced function $g_q$ has exactly one cycle $C$, and let
$I=(I_1,\cdots,I_n)$ be an initial word such that
$(I_1,\cdots,I_{n-1})$ is on $C$. Then $(P,I)$ is a de~Bruijn
sequence.
\end{corollary}

\section{Generalizing the prefer-Opposite and prefer-same sequences}\label{S:span1}
In this section we will focus on preference functions $P$ of span
$1$, that can be represented by $q\times q$ a matrices, where for
each digit $j=0,\ldots,q-1$, $P(j)$ is given by the $(j+1)^{st}$ row
of the matrix. When $q=2$, the only possible matrices are

\[
F_2=\begin{bmatrix} 1 & 0\\1 & 0\end{bmatrix}; O_2=\begin{bmatrix} 1 & 0\\0 &
1\end{bmatrix}; S_2=\begin{bmatrix} 0 & 1\\1 & 0\end{bmatrix};\;\; \text{
and } Z_2=\begin{bmatrix} 0 & 1\\0 & 1\end{bmatrix}.
\]
Applying Theorem~\ref{T:main} to the first matrix with initial word
$0^n$, the set $\Sigma$ is empty and we obtain the famous prefer-one
sequence, first attributed to Martin~\cite{Martin1934}. With the
same initial word, the second matrix gives $q^{\prime}=2$ and
$\Sigma=\{1\}$, leading to the prefer-opposite sequence, see
Alhakim~\cite{Alhakim10}. For the third sequence, starting with one
of the alternating strings $010\cdots$ or $101\cdots$ leads to a
full sequence. \cite{Alhakim19} shows that these sequences give a
cyclic rotation of the prefer-same sequence (and bitwise
complement), thus presenting a new and simplified way to generate
the latter sequence, avoiding some subtleties in the original
algorithm described in Eldert~\cite{Eldert} and
Fredricksen~\cite{Fred1982}. The fourth matrix, the `prefer-zero'
matrix, evidently produces the complement of the prefer-one sequence
when the initial word is $1^n$.

As mentioned in the introduction, the prefer-opposite is already
generalized to the so-called prefer-higher sequence for general
alphabet size  $q$, see Fredricksen~\cite{Fred1982}. For the other
two sequences, the concept of `opposite' is of course not naturally
meaningful. However, we will employ preference functions of span
$1$--representable by $q\times q$ matrices--to obtain natural
generalization of the prefer-opposite and prefer-same sequences to
the case of alphabet size $q$.

Fix a digit $d\in \mathcal{A}$. We define the two $q\times q$
matrices $O:=O^{(q,d)}$ and $S:=S^{(q,d)}$ by

$O_{ij}=i+(j+1)d\mod q$ for $i,j=0,1,\ldots,q$, and

$S_{ij}=i+(j)d\mod q$ for $i,j=0,1,\ldots,q$

Examples of preference functions of the $O$ type are given in the
top row of Table~\ref{Ta:GenPref_Diagrams}, while the lower row
illustrates preference functions of Type $S$.

\begin{lemma}\label{L:co-prime}
If $gcd(d,q)=1$, then each of the matrices $O$ and $S$ defines a
preference function of span 1; the column functions for both
preference functions are bijective functions from $\mathcal{A}$ to
$\mathcal{A}$. Moreover, for all digits $a\in \mathcal{A}$,
$O_q(a)=a$  and $S_1(a)=a$.
\end{lemma}

\begin{proof}
For each $i$, the $i^{th}$ row of $O$ is $(i+d,i+2d,\cdots,i+qd)$,
and the $i^{th}$ row of $S$ is $(i,i+d,i+2d,\cdots,i+(q-1)d)$; where
addition modulo $q$. Since $gcd(d,q)=1$, it is immediate that the
$q$ entries are distinct in both scenarios, verifying that both $O$
and $S$ are preference functions of span 1. Also in both scenarios,
every column function has the form $g_j(i)=i+jd\mod q$ for some
fixed value of $j=1,\ldots,q$, which is evidently bijective. The
other two statements of the lemma are obvious.
\end{proof}

\begin{theorem}~\label{T:pref_opp}
The sequence $\textbf{o}_n=(O,0^n)$, with $gcd(d,q)=1$, is a
de~Bruijn sequence missing only the words $1^n,\ldots,(q-1)^n$.
\end{theorem}

\begin{proof}
We apply Theorem~\ref{T:main} with $C$ being the self loop
$0\rightarrow0$. By Lemma~\ref{L:co-prime}, $O_q(a)=a$ for all $a\in
\mathcal{A}$, so that $\Sigma=\mathcal{A}\backslash\{0\}$ is closed
under $g_q$. Since $O_{q}(0)=0$, $O_{q-1}(0)\neq0$. $g_{q-1}$ is a
bijection (by Lemma~\ref{L:co-prime}) so there exists $b\in\Sigma$
such that $g_{q-1}(b)=0$, which is not in $\Sigma$. Therefore,
$q^{\prime}:=q$ is minimal in the sense of Theorem~\ref{T:main}. We
now verify that $g_{q^{\prime}-1}$ has no cycles when restricted to
$\Sigma$. Noting that $gcd(q,d(q-1))=1$,
$g_{q^{\prime}-1}(i)=i+(q-1)d$ induces exactly one cycle. Thus,
$g^q_{q-1}(0)=0$. Let $c\in\Sigma$ be such $g_{q-1}(0)=c$. Iterating
$g_{q-1}$ we get

$(c,g_{q-1}(c), g_{q-1}^2(c),\cdots,g_{q-1}^{q-1}(c))$

which forms a ``path'' whose elements are all different, exhausting
$\Sigma$ and with $0$ at the end. This verifies that $g_{q-1}$ has
no cycles in $\Sigma$. By Theorem~\ref{T:main}, the only missing
words are as claimed.
\end{proof}

\begin{definition}\label{D:Alternating}
We say that $(y_1,\ldots,y_n)$ is an alternating string with step
$l$, $l=1,\ldots,q-1$ if $y_{i+1}=y_i+ld \mod q$ for
$i=1,\ldots,n-1$.
\end{definition}

\begin{theorem}~\label{T:pref_same}
Let $I$ be an alternating string with step $d(q-1)$ and
$gcd(d,q)=1$. Then the sequence $\textbf{s}_n=(S,I)$  is a full
de~Bruijn sequence.
\end{theorem}

\begin{proof}
Note that $S_q(a)=a+d(q-1)$ for any $a\in \mathcal{A}$. Since
$d(q-1)$ is relatively prime to $q$, the induced function $g_q$
cycles through all symbols in $\mathcal{A}$, i.e. words of size $s=1$.
It follows that the initial $I$ is on the cycle of $g_q$. Thus the
statement follows immediately from Corollary~\ref{C:complete}.
\end{proof}

The following gives an interpretation of the sequence $\textbf{o}_n$
as a prefer-opposite sequence. If we \emph{define} the opposite of
$i$ to be $i+d\mod q$, then in the construction of the sequence we
first propose the opposite of the last produced digit. If the
resulting $n$-word is new we accept the proposed digit. Otherwise we
propose the opposite of the last rejected digit and so. We stop when
$q$ proposals for the same digit are rejected. Noting that $S$ can
be obtained from $O$ by moving the last column of the latter to the
first column, the preference function defined by $S$ gives the
highest preference to the most recent digit (hence prefer-same) and
then proceeds by skipping the value $d$ each time a proposed digit
is rejected.

Two notions of \textit{opposite} of $i$ which are available for any
$q$ are the next digit $i+1$, and the \textit{farthest} digit
$i+q-1$. For $q=5$, these are illustrated in the top part of
Table~\ref{Ta:GenPref_Diagrams},
along with a third case, $d=2$, in the middle. 

\begin{table}[h!]
\small
\begin{tabular}{c|c|c}
\hline
\begin{tabular}{lllllll}
0 & $\rightarrow$ 1, & 2, & 3, & 4, & 0\\
1 & $\rightarrow$ 2, & 3, & 4, & 0, & 1\\
2 & $\rightarrow$ 3, & 4, & 0, & 1, & 2\\
3 & $\rightarrow$ 4, & 0, & 1, & 2, & 3\\
4 & $\rightarrow$ 0, & 1, & 2, & 3, & 4\\
\end{tabular}
&
\begin{tabular}{lllllll}
0 & $\rightarrow$ 2, & 4, & 1, & 3, & 0\\
1 & $\rightarrow$ 3, & 0, & 2, & 4, & 1\\
2 & $\rightarrow$ 4, & 1, & 3, & 0, & 2\\
3 & $\rightarrow$ 0, & 2, & 4, & 1, & 3\\
4 & $\rightarrow$ 1, & 3, & 0, & 2, & 4\\
\end{tabular}
&
\begin{tabular}{lllllll}
0 & $\rightarrow$ 4, & 3, & 2, & 1, & 0\\
1 & $\rightarrow$ 0, & 4, & 3, & 2, & 1\\
2 & $\rightarrow$ 1, & 0, & 4, & 3, & 2\\
3 & $\rightarrow$ 2, & 1, & 0, & 4, & 3\\
4 & $\rightarrow$ 3, & 2, & 1, & 0, & 4\\
\end{tabular}\\\hline\hline
\begin{tabular}{lllllll}
0 & $\rightarrow$ 0, & 1, & 2, & 3, & 4\\
1 & $\rightarrow$ 1, & 2, & 3, & 4, & 0\\
2 & $\rightarrow$ 2, & 3, & 4, & 0, & 1\\
3 & $\rightarrow$ 3, & 4, & 0, & 1, & 2\\
4 & $\rightarrow$ 4, & 0, & 1, & 2, & 3\\
\end{tabular}
&
\begin{tabular}{lllllll}
0 & $\rightarrow$ 0, & 2, & 4, & 1, & 3\\
1 & $\rightarrow$ 1, & 3, & 0, & 2, & 4\\
2 & $\rightarrow$ 2, & 4, & 1, & 3, & 0\\
3 & $\rightarrow$ 3, & 0, & 2, & 4, & 1\\
4 & $\rightarrow$ 4, & 1, & 3, & 0, & 2\\
\end{tabular}
&
\begin{tabular}{lllllll}
0 & $\rightarrow$ 0, & 4, & 3, & 2, & 1\\
1 & $\rightarrow$ 1, & 0, & 4, & 3, & 2\\
2 & $\rightarrow$ 2, & 1, & 0, & 4, & 3\\
3 & $\rightarrow$ 3, & 2, & 1, & 0, & 4\\
4 & $\rightarrow$ 4, & 3, & 2, & 1, & 0\\
\end{tabular}
\\\hline
\end{tabular}
\caption{Top row: generalized prefer-opposite diagrams with
$d=1,2,4$. Lower row: generalized prefer-same diagrams with
$d=1,2,4$.}
\label{Ta:GenPref_Diagrams}
\end{table}

Besides this intuitive interpretation, the proposed sequences bear
more resemblance to the binary counterparts, as we now present. We
first focus on $\textbf{o}_n$ and later present similar results for
$\textbf{s}_n$. We recall here that the binary prefer-opposite
sequence of order $n$ has all words except $1^n$, and it terminates
with the substring $1^{n-1}0^{n-1}$, see Alhakim~\cite{Alhakim10}.
Proposition~\ref{P:opp_suffix} below shows how this phenomenon
generalizes in the current setting.

\begin{lemma}~\label{L:FinalAppearance}
In the sequence $\textbf{o}_n$, the digit $a$ does not appear after
the last occurrence of the pattern $a^{n-1}$ for any  $a\in
\mathcal{A}$.
\end{lemma}

\begin{proof}
The statement is trivial when $a=0$ because $0^{n-1}$ is the
terminal pattern. Suppose the statement is false for some $a\neq0$
and suppose that $x_1\cdots x_{n-1}a$ is the first $n$-word that
terminates with $a$ after the last occurrence of $a^{n-1}$.

It follows that $x_1\cdots x_{n-1}a$ must be followed by something,
as it is not the last word. The next $n$-word is either $x_2\cdots
x_{n-1}aa$, or otherwise the word $x_2\cdots x_{n-1}aa$ will occur
afterwards in the sequence, because $O_q(a)=a$.

The same argument can now be iterated to establish that $a^{n-1}$
will appear later in the sequence, thus contradicting the assumption
that $a^{n-1}$ appeared for the last time earlier.
\end{proof}

\begin{proposition}\label{P:opp_suffix}
For any integer $d$ co-prime with $q$, the generalized
prefer-opposite sequence $\textbf{o}_n$ terminates with the pattern
$$\sigma_{q-2}^{n-1}\cdots\sigma_{1}^{n-1}\sigma_{0}^{n-1}0^{n-1}$$
where $\sigma_0,\cdots,\sigma_{q-2}$ are an appropriate
re-arrangement of the digits $1,\cdots q-1$ determined by the
function $g_{q-1}$.
\end{proposition}

\begin{proof}
Since $g_{q-1}$ is cyclic, the terms of the sequence
\[
g_{q-1}(0),g_{q-1}^2(0),\ldots,g_{q-1}^q(0)=0
\]
are all distinct, where $g^{i+1}_{q-1}(0)=g^{i}_{q-1}(g_{q-1}(0))$.
Let $a=g^{i}_{q-1}(0)$ and $b=g^{i+1}_{q-1}(0)=g_{q-1}(a)$ for an arbitrary $i<q-1$, so that $a$ and $b$ are distinct from $0$. By the preference order,
the last occurrence of $a^{n-1}$ is succeeded by the digit $b$. We
will show that it is in fact succeeded by $b^{n-1}$.

Let $j$ be the largest integer such that $a^{n-1}$ is succeeded by
$b^j$. By Theorem~\ref{T:pref_opp}, $1\leq j\leq n-1$ so suppose
$j<n-1$. It then follows that the $n$-word $a^ib^j$ that is a suffix
of $a^{n-1}b^j$ has $i\geq2$. Since $a^ib^j$ is not the terminal
$n$-word of the sequence it must be succeeded by a symbol $c$. Since
$j$ is maximal, $c\neq b$. This puts the word $a^{i-1}b^jc$ and it
follows that $a^{i-1}b^{j+1}$ will appear later in the sequence
because $O_q(b)=b$. This contradicts Lemma~\ref{L:FinalAppearance}
as $i-1\geq1$, therefore $j=n-1$.

Letting $\sigma_{q-1-i}=g_{q-1}^i(0)$, the previous paragraph
establishes that the last occurrence of $\sigma_{q-1-i}^{n-1}$ is
succeeded by $\sigma_{q-1-i-1}^{n-1}$ for $i=1,\ldots,q-2$. When
$i=q-1$, observe that $g_{q-1}(\sigma_0)=0$, so the last occurrence
of $\sigma_0^{n-1}$ is followed by $0$. At this point no nonzero
symbol  is allowed to appear, by Lemma~\ref{L:FinalAppearance}. Thus
the only remaining possibility is the terminal string $0^{n-1}$ and
this finishes the proof.
\end{proof}

As an example, the order $n$ sequences generated by the generalized
prefer-opposite diagrams in Table~\ref{Ta:GenPref_Diagrams}
terminate, respectively from left to right, with the patterns

$4^{n-1}3^{n-1}2^{n-1}1^{n-1}0^{n-1}$,
$3^{n-1}1^{n-1}4^{n-1}2^{n-1}0^{n-1}$, and
$1^{n-1}2^{n-1}3^{n-1}4^{n-1}0^{n-1}$.

\begin{remark}
The case $n=2$ yields the following interesting palindromic sequences\\

$0012340241303142043210$,\hspace{5pt} $0024130432101234031420$,
\hspace{5pt}
$0043210314202413012340$.\\

The reader is invited to check that a palindrome is obtained (with
another zero appended at the end) if and only if $n=2$ and $q$ is a
prime number.
\end{remark}

\subsection{Mapping $\textbf{o}_n$ to the prefer-higher sequence}

The most striking similarity with the binary prefer-opposite is the
relationship between the proposed sequence and the prefer-higher
sequence, which is essentially the same relationship between the
binary prefer-opposite and the prefer-one sequences, discussed
by~\cite{RW2017}.

To describe this relationship let us first recall the so-called
$D$-homomorphism used in Lempel~\cite{Lempel1970} to recursively
construct binary de~Bruijn sequences. For two binary numbers $b_1$
and $b_2$ it is defined as $D(b_1,b_2)=b_1+b_2\mod 2$. It can more
generally be defined as a function from $\{0,1\}^l$ to
$\{0,1\}^{l-1}$, for any integer $l\geq2$ as follows,
$D(b_1,\ldots,b_l)=(b_1+b_2,b_2+b_3,\ldots,b_{l-1}+b_l)$, where
addition is performed modulo 2. In fact, Rubin and
Weiss~\cite{RW2017} establish that if we map the binary
prefer-opposite sequence $(b_0,\ldots, b_{2^n-1})$ via the
$D$-homomorphism, getting $(\hat{b}_0,\ldots, \hat{b}_{2^n-2})$, and
then drop all bits $\hat{b}_i$ such that
$(\hat{b}_{i-n+2},\ldots,\hat{b}_{i})$ is already a substring of
$(\hat{b}_0,\ldots,\hat{b}_{i-1})$, the remaining bits form the
prefer-one sequence of order $n-1$.

To state the $q$-ary counterpart of this result, we will use the
generalized $q$-ary $D$-homomorphisms defined in Alhakim and
Akinwande~\cite{AA2011} as

$$D_{\beta}(x_1,\ldots,x_l)=\left(\beta(x_2-x_1),\ldots,\beta(x_l-x_{l-1})\right),$$

where $\beta$ is co-prime with $q$, subtraction is performed modulo
$q$, and $x_i$ is a digit in the $q$-ary alphabet $\mathcal{A}$.
These homomorphisms were also used to generalize the Lempel
construction to nonbinary alphabets, see~\cite{AA2011}.

\begin{theorem}~\label{T:MappingOpp2Hi}
Let
$D_{\beta}\left(\bf{o}_n\right)=(\hat{x}_0,\ldots,\hat{x}_{q^n-2})$
be the homomorphic image of the generalized prefer-opposite sequence
${\bf{o}_n}=(x_0,\ldots,x_{q^n-1})$ by the homomorphism $D_{\beta}$,
where $gcd(d,q)=1$ and $\beta=d^{-1}(q-1)$, with $d^{-1}$ being the
inverse of $d$ modulo $q$.

\noindent Then the subsequence of $D_{\beta}\left(\bf{o}_n\right)$
obtained by deleting every digit $\hat{x}_i$ such that
$(\hat{x}_{i-n+2},\ldots,\hat{x}_{i})$ is a substring of
$(\hat{x}_0,\ldots,\hat{x}_{i-1})$, is the prefer-higher de~Bruijn
sequence ${\bf{h}_{n-1}}$ of order $n-1$.
\end{theorem}

To help the reader follow the proof more closely, we illustrate this
result with an example. For $q=3$, $d=2$, $\bf{o}_4$ is

\noindent{\small
$0000210212101020211022100201012120200220222122112111011001000120122011200111222000$}

As $\beta=1$, the appropriate homomorphism is
$D_{1}(x_1,x_2)=x_2-x_1$, giving the image sequence
$D_1(\bf{o}_4)$:\newline

\noindent{\small
$000222221221221220220220211211211210201200210201200210201200111110110110100100100$}\newline

To best present how the prefer-higher sequence of order $3$ emerges
from the latter sequence, we write it as a sequence of words and
cross out any repetition, keeping only the first occurrence of every
word.

$000\rightarrow002\rightarrow022\rightarrow222\rightarrow\cancel{222}\rightarrow\cancel{222}%
\rightarrow221\rightarrow212\rightarrow122\rightarrow\cancel{221}\rightarrow\cancel{212}%
\rightarrow\cancel{122}\rightarrow\cancel{221}\rightarrow\cancel{212}\rightarrow\cancel{122}%
\rightarrow220\rightarrow202\rightarrow\cancel{022}\rightarrow\cancel{220}\rightarrow\cancel{202}%
\rightarrow\cancel{022}\rightarrow\cancel{220}\rightarrow\cancel{202}\rightarrow021\rightarrow211%
\rightarrow112\rightarrow121\rightarrow\cancel{211}\rightarrow\cancel{112}\rightarrow\cancel{121}%
\rightarrow\cancel{211}\rightarrow\cancel{112}\rightarrow\cancel{121}\rightarrow210\rightarrow102%
\rightarrow020\rightarrow201\rightarrow012\rightarrow120\rightarrow200\rightarrow\cancel{002}%
\rightarrow\cancel{021}\rightarrow\cancel{210}\rightarrow\cancel{102}\rightarrow\cancel{020}%
\rightarrow\cancel{201}\rightarrow\cancel{012}\rightarrow\cancel{120}\rightarrow\cancel{200}%
\rightarrow\cancel{002}\rightarrow\cancel{210}\rightarrow\cancel{102}\rightarrow\cancel{020}%
\rightarrow\cancel{201}\rightarrow\cancel{012}\rightarrow\cancel{120}\rightarrow\cancel{200}%
\rightarrow001\rightarrow011\rightarrow111\rightarrow\cancel{111}\rightarrow\cancel{111}%
\rightarrow110\rightarrow101\rightarrow\cancel{011}\rightarrow\cancel{110}\rightarrow\cancel{101}%
\rightarrow\cancel{011}\rightarrow\cancel{110}\rightarrow\cancel{101}\rightarrow010\rightarrow100%
\rightarrow\cancel{001}\rightarrow\cancel{010}\rightarrow\cancel{100}\rightarrow\cancel{001}%
\rightarrow\cancel{010}\rightarrow\cancel{100}$

\noindent The uncrossed words can be written compactly as the
sequence $$00022212202112102012001110100,$$ which is the prefer
higher sequence of order $3$.


\begin{definition}
We say that $(y_1,\ldots,y_n)$ is a translate of $(x_1,\ldots,x_n)$
if there exists  $c\in \mathcal{A}$ such that $y_i=x_i+c \mod q$ for
$i=1,\ldots,n$.
\end{definition}

\begin{lemma}\label{L:translates}
Given two $n$-tuples $\bf{v}$ and $\bf{w}$,
$D_{\beta}({\bf{v}})=D_{\beta}(\bf{w})$ if and only if $\bf{w}$ is a
translate of $\bf{v}$.
\end{lemma}

\begin{lemma}\label{L:q-inverses}
Every $(n-1)$-tuple has exactly $q$ inverse images by $D_{\beta}$.
\end{lemma}

\begin{proof}
Let ${\bf{w}}_{n-1}=(\hat{x}_1,\ldots,\hat{x}_{n-1})$. For each
$x_1$ in $\mathcal{A}$, define
$x_i=x_{i-1}+\beta^{-1}\hat{x}_{i-1}$, $i=2,\ldots,n$. One can
readily check that $D_{\beta}(x_1,\ldots,x_n)={\bf{w}}_{n-1}$.
\end{proof}

In general, there exist $q/gcd(q,l)$ alternating strings with step
$l$; see Definition~\ref{D:Alternating}, which are all translates of
each other. Observe that the $q$ alternating strings with step $l=1$
occur consecutively in the beginning of the sequence. Interestingly,
alternating strings of higher steps also aggregate together.

\begin{lemma}\label{L:Alternating}
When the first alternating string with step $l;\;l>1$ occurs,
all of its translates follow consecutively.
\end{lemma}

\begin{proof}
Let $l=\hat{l}d$, where $\hat{l}$ is between $1$ and $q-1$,
inclusive. Due to the form of the preference matrix $O^{(q,d)}$, an
alternating string with step $l$ corresponds to a cycle of
$g_{\hat{l}}$. That is, when the first digit of an alternating
string with step $l$ is given, all other digits are entered via the
$\hat{l}^{th}$
preference (the $1^{st}$ preference being the highest). By the prefer-%
opposite strategy, the first alternating string is achieved when all
words that have combinations of preference steps less than or equal
to $\hat{l}$, but not purely $\hat{l}$, have been used. Immediately
after this, preference step $\hat{l}$ will be repeated until a digit
of preference $\hat{l}$ cannot be added without putting an already
present word. This is the case when a cycle of the function
$g_{\hat{l}}$ is exhausted.

There are two cases. If $\hat{l}$ is co-prime with $q$, then
$g_{\hat{l}}$ is periodic with period $q$, in which case all
alternating strings of step $l$ are placed consecutively. Otherwise,
preference $\hat{l}+1$ is used once and then $\hat{l}$ is used
again, that is, another cycle of $g_{\hat{l}}$ is entered, and so on
until all cycles of $g_{\hat{l}}$ are placed, thus producing all
alternating strings with step $l$.
\end{proof}

\begin{lemma}\label{L:Opp2Hi}
A word $\textup{\bf{w}}=(x_1,\ldots,x_n)$ is succeeded by the digit
$x_n+jd$ for some $j=0,1,\ldots,q-1$ and $d$ that is co-prime with
$q$, \textit{if and only if} $D_{\beta}(x_1,\cdots,x_n)$ is
succeeded by $q-j\mod q$.
\end{lemma}

\begin{proof}
$D_{\beta}\left(x_1,\cdots,x_n\right)=\left(\beta(x_2-x_1),\ldots,\beta(x_n-x_{n-1})\right)$.
If $x_n+jd$ succeeds $\bf{w}$, $D_{\beta}(\textbf{w})$ is succeeded
by digit $\beta(x_n+jd-x_n)=\beta jd$. Using the value of $\beta$ we
get $d^{-1}(q-1)jd=j(q-1)=q-j$.

Conversely, suppose $D_{\beta}(x_1,\cdots,x_n)$ is succeeded by
$q-j\mod q$ and let $y$ be the digit that succeeds
$\textup{\bf{w}}$. That is, $\beta(y-x_n)=q-j$. This implies that
$y-x_n=\beta^{-1}(q-j)=d(q-1)(q-j)=dj\mod q$ and the result follows.
\end{proof}

Lemma~\ref{L:Opp2Hi} provides the intuition behind the result of
Theorem~\ref{T:MappingOpp2Hi}. In the prefer-opposite sequence, the
$q$ appearances of an $(n-1)$-tuple $(x_1,\ldots,x_{n-1})$ in
$\textbf{o}_n$ are succeeded respectively by
$x_{n-1}+d,x_{n-1}+2d,\ldots,x_{n-1}+qd=x_{n-1}$. Their images by
$D_\beta$ are followed respectively by $q-1,\ldots,1,0$, i.e.,
according to the prefer-higher strategy. However, the problem is
more complicated. Using Lemma~\ref{L:translates}, we see that a
translate of the former $n$-tuples can potentially occur in
$\textbf{o}_n$ in the wrong location so as to mess up the
prefer-higher order.

Furthermore, another main issue is to verify that, upon crossing out
a sequence of one or more consecutive words, the last word before
this crossed sequence teams up correctly with the next uncrossed
word.

In the binary case, Rubin and Weiss~\cite{RW2017} carry out a
meticulous study to inspect various cases of which translates of a
word occur (in this case the only translate of a word is its bitwise
complement.)  Our proof will apply the principle of strong
mathematical induction, which greatly simplifies the discussion.


\begin{proof}[Proof of Theorem \ref{T:MappingOpp2Hi}]
Observe that the inverse images of ${\bf{0}_{n-1}}$ are all words of
the form $c^n$. By Theorem~\ref{T:pref_opp}, $0^n$ is the only
inverse image of $0^{n-1}$ that occurs in ${\bf{o}_n}$. Moreover,
any word that is not of the form $c^n$ occurs in $\bf{o}_n$ along
with all of its translates.

By Lemmas~\ref{L:translates} and \ref{L:q-inverses}, the collection
of inverse images of all $n-1$ tuples is the collection of $n$
tuples. It follows that $D_{\beta}\left(\bf{o}_n\right)$ includes
every $n-1$ tuple exactly $q$ times, except for $0^{n-1}$ that
appears once in the beginning.

Hence, there will be exactly $q^{n-1}$ uncrossed digits. We let
$I=\left(i_1,i_2,\cdots,i_{q^{n-1}}\right)$ be the subsequence of
$1,\ldots,q^n-1$ consisting of the indices in ${\bf{\hat{o}}_{n-1}}$
that remain after crossing repetitions. That is, $i_t\in I$ if and
only if
$\left(\hat{x}_{i_t},\ldots,\hat{x}_{i_t+n-2}\right)\notin\{(\hat{x}_j,\ldots,\hat{x}_{j+n-2}):j=1,\ldots,i_t-n+1\}$.

Let us also represent the prefer-higher sequence ${\bf{h}_{n-1}}$ as
$\left(y_i;i=0,\ldots,q^{n-1}-1\right)$, and define
$\textbf{x}_{i}:=\left(x_{i},\ldots,x_{i+n-1}\right)$,
$\hat{\textbf{x}}_{i}:=\left(\hat{x}_{i},\ldots,\hat{x}_{i+n-2}\right)$
and $\textbf{y}_{i}:=\left(y_{j},\ldots,y_{j+n-2}\right)$ for
$i\geq0$ to be the $i^{th}$ word of $\bf{o}_n$,
$D_{\beta}\left(\bf{o}_n\right)$ and $\bf{h}_{n-1}$ respectively.
Note that $D_{\beta}(\textbf{x}_{i})=\hat{\textbf{x}}_{i}$.

We need to verify that (1) for each $t$ the last $n-2$ digits
$\hat{\textbf{x}}_{i_t}$ coincide with the first $n-2$ digits of
$\hat{\textbf{x}}_{i_{t+1}}$, and (2) the prefer higher strategy is
respected. A moment's reflection shows that it is sufficient to show
that $\hat{\textbf{x}}_{i_t}=\textbf{y}_{t}$ for all $t;\;0\leq
t\leq q^{n-1}-1$. Observing that this statement is true for $t=0$,
we assume that it is true for all $j\leq t$, where $0\leq t<
q^{n-1}-1$ and proceed to prove it is true for $t+1$.

First consider the case when $i_{t+1}=i_{t}+1$. Suppose that
$\hat{\textbf{x}}_{i_t+1}$ ends with the digit $q-j$ for some $j$,
then by Lemma~\ref{L:Opp2Hi} $\textbf{x}_{i_t+1}$ has the form
$(x_{i_t+1},\ldots,x_{i_t+n-1},x_{i_t+n-1}+jd)$. By the
prefer-opposite strategy, the $j-1$ words
$(x_{i_t+1},\ldots,x_{i_t+n-1},x_{i_t+n-1}+ld)$, $1\leq l<j$ occur
earlier in ${\bf{o}_n}$. Since they all agree with
$\textbf{x}_{i_t+1}$ on the first $n-1$ digits, their images by
$D_{\beta}$ are of the form
$(\hat{x}_{i_{t}+1},\ldots,\hat{x}_{i_{t}+1+n-2},q-l)$ and they all
occur earlier. To see that $\hat{\textbf{x}}_{i_t+1}$ will not be
crossed out, suppose that it did occur as $\textbf{y}_s,\;s\leq t$.
So its first occurring inverse image in ${\bf{o}_n}$ must be
$\textbf{x}_{i_s}$, by the induction hypothesis, so that
$\textbf{x}_{i_t+1}$ is not new, contradicting the fact that
$i_t+1=i_{t+1}$.

When $i_{t+1}>i_{t}+1$, $i_{t+1}-1\notin I$ and
$D_{\beta}(\textbf{x}_{i_{t+1}-1})=\textbf{y}_{\tau}$ for some
$\tau\leq t$. By the inductive hypothesis this means that
$\textbf{y}_{\tau}=D_{\beta}(\textbf{x}_{i_{\tau}})$. There are $q$
possible $(n-1)$-tuples that can logically succeed
$\textbf{y}_{\tau}$:
$\textbf{w}_l=(y_{\tau+1},\ldots,y_{\tau+n-2},q-k);\;k=1,\ldots,q$.

Since $i_{\tau}\leq i_t$, 
$i_{\tau}+1\leq i_{t+1}-1$ so that $i_{\tau+1}\notin I$. Thus at
least one of the $\bf{w}_l$ occurs in
$\hat{\textbf{x}}_{i_0},\ldots,\hat{\textbf{x}}_{i_t}$. Let $j$,
$1\leq j\leq q$ be the maximum integer such that $\textbf{w}_j$
occurs. Then by the inductive hypothesis, all
$\textbf{w}_l,\;\emph{l}\leq j$ occur before since so far
$\hat{\textbf{x}}_{i_0},\ldots,\hat{\textbf{x}}_{i_t}$ coincide with
$\textbf{y}_0,\ldots,\textbf{y}_t$. This also implies that $j<q$,
for otherwise $i_{t+1}\notin I$ (as
$D_{\beta}(\textbf{x}_{i_{t+1}})$ has to go to one of the
$\textbf{w}_l$ and has to be new.) It also follows that
$D_{\beta}(\textbf{x}_{i_{t+1}})$ is $\textbf{w}_{j+1}$, for if it
equals $\bf{w}_k$ for some $k>j+1$, $\textbf{x}_{i_{t+1}}$ would
have the form
$(x_{i_{t+1}},\ldots,x_{i_{t+1}+n-2},x_{i_{t+1}+n-2}+kd)$. The
prefer-opposite strategy would then mean that
$(x_{i_{t+1}},\ldots,x_{i_{t+1}+n-2},x_{i_{t+1}+n-2}+(j+1)d)$ occurs
earlier in ${\bf{o}_n}$ and therefore $\textbf{w}_{j+1}$ occurs
earlier in $D_{\beta}({\bf{o}_n})$, contradicting the maximality of
$j$.

To sum up, the prefer-opposite again implies that the $n$-tuples

\begin{equation}~\label{E:Xkappa}
\textbf{x}_{\kappa_l}=(x_{i_{t+1}},\ldots,x_{i_{t+1}+n-2},x_{i_{t+1}+n-2}+ld),\;1\leq
l\leq j
\end{equation}

all occur earlier in ${\bf{o}_n}$, and
$\{D_{\beta}(\textbf{x}_{\kappa_l}):1\leq l\leq
j\}=\{\textbf{w}_1,\ldots,\textbf{w}_j\}$. The
$\textbf{x}_{\kappa_l}$ are evidently all distinct, and there are
two possibilities.

Case 1. Suppose $\textbf{x}_{i_{\tau}}\neq \textbf{x}_{\kappa_l}$
for any $l=1.\ldots,j$. If $\kappa_1,\ldots,\kappa_j, i_{\tau}$ are
all strictly less than $i_t$ then their images by $D_{\beta}$ are
$\textbf{w}_1,\ldots,\textbf{w}_{j+1}$, contradicting that
$i_{t+1}\in I$. Thus $i_t\in\{\kappa_1,\ldots,\kappa_j,i_{\tau}\}$,
and by Equation~(\ref{E:Xkappa}) and the prefer-opposite,
$i_{\tau}=i_t$.

Case 2. If $\bf{x}_{i_{\tau}}= \bf{x}_{\kappa_l}$ for some $l$, then
$\bf{x}_{i_{\tau}}$ has the form of an alternating string with step
$l$. Since $\tau\in I$, $\bf{x}_{i_{\tau}}$ is the first such
alternating string. By Lemma~\ref{L:Alternating}, all alternating
strings of step $l$ and that are translates of each other  occur
consecutively. This forces the last of these strings to be
$\bf{x}_{i_{t+1}-1}$ and thus forces $i_{\tau}$ to be $i_t$, i.e.
$t=\tau$.

\end{proof}

\section{Discrepancy Calculations}
Cooper and Heitsch \cite{Cooper2010} define the discrepancy of a
binary sequence $w=w_1\cdots w_n$ as
$disc(w)=\underset{M}{max}|\displaystyle\sum_{k=0}^{M}(-1)^{w_k}|$,
where the maximum is taken over all initial strings of length $M$.
In fact, they prove that the discrepancy of the prefer-one sequence
of word size $n$ is of order $\Theta(2^n\ln n/n)$, much larger than
that of a typical de~Bruijn sequence. This is the case because the
prefer-one scheme biases the initial part of the sequence (started
with $0^n$) to having more ones than zeros.

Using a variant of the above discrepancy formula, Gabric and
Sawada~\cite{GS2022} numerically studied the discrepancy of a number
of binary de~Bruijn sequences, including the prefer-same and
prefer-opposite sequences which come out to have a low discrepancy.
In this section we define a measure of discrepancy for sequences of
general alphabet size and we use it to compare the prefer-higher
sequence and the two proposed sequences with the binary versions.

\begin{definition}
Given is a finite sequence $s$ from an alphabet $\mathcal{A}$ of
size $q$. For a letter $a\in \mathcal{A}$, we write $f_a(i)$ to
denote the frequency of $a$ in the initial subsequence of size $i$.
The discrepancy of $s$ is defined as $$discrep(s)=
\underset{i}{max}\left(\underset{a,b}{max}|f_a(i)-f_b(i)|\right),$$
where the inner maximum is taken over all pairs of digits $a$ and
$b$, while the outer maximum is taken over all initial subsequences
of size $i$ of the sequence.
\end{definition}


We observe that this discrepancy, when applied to a de Bruijn
sequence (or any periodic sequence), is a function of the initial
word. When $q=2$, this definition boils down to the formula given in
\cite{Cooper2010}.

Gabric and Sawada~\cite{GS2022} conjecture that the prefer-opposite
and prefer-same sequences have a discrepancy in the order of
$\Theta(n^2)$. Table~3 below displays the calculated discrepancy for
all three sequences for $q=2,\ldots,5$ and for several values of
$n$, when each sequence starts with its natural, cyclic initial
word. Notably, both proposed sequences have discrepancies that are
very similar to their binary counterparts, regardless of the value
of $q$ and $n$. Also, the table strongly suggests that the
discrepancy of both sequences are dominated by $(n-1)^2$ while that
of the prefer-higher sequence grows essentially exponentially, for
various values of $q$. One can even guess the patterns for the first
two sequences. When $q=2$ for instance, the discrepancy of the
pref-same, for $n$ high enough, is
\[discrep(\text{pref-same}(n))=\begin{cases} n(n-2)/4 & n\text{ is
even}\\(n-1)^2/4 & (n-1)\mod 4=0\\(n-1)^2/4+1 & (n-1)\mod
4=2\\\end{cases},\] while that of the prefer-opposite appears to be
\[discrep(\text{pref-opp}(n))=\begin{cases} (\frac n2)^2+1 & n\text{ is
even}\\(\frac{n-1}2)^2+\frac{n-1}2+1 & n\text{ is
odd}\\\end{cases}.\] Curiously also,
$discrep(\text{pref-opp}(n))=(n-1)^2$ for $n$ high enough when
$q=4,5$.

\begin{table}
\center
\begin{tabular}{l||l|l|l||l|l|l||l|l|l||l|l|l}
\multicolumn{13}{c}{Discrepancy(pref-same, pref-opp,
pref-higher)}\\\hline $n$ & \multicolumn{3}{c||}{$q=2$} &
\multicolumn{3}{c||}{$q=3$} & \multicolumn{3}{c||}{$q=4$} &
\multicolumn{3}{c}{$q=5$}\\\hline
2 & 1 & 2 & 2 & 2 & 2 & 2 & 3 & 3 & 3 & 3 & 2 & 4\\
3 & 3 & 3 & 3 & 3 & 4 & 5 & 4 & 5 & 9 & 6 & 4 & 16\\
4 & 3 & 5 & 4 & 5 & 7 & 12 & 6 & 9 & 31 & 9 & 9 & 70\\
5 & 5 & 7 & 5 & 8 & 12 & 30 & 9 & 16 & 112 & 11 & 16 & 310\\
6 & 6 & 10 & 6 & 10 & 19 & 84 & 13 & 25 & 391 & 14 & 25 & 1346\\
7 & 10 & 13 & 14 & 16 & 27 & 227 & 18 & 36 & 1407 & 17 & 36 & 5807\\
8 & 12 & 17 & 28 & 21 & 37 & 618 & 23 & 49 & 5152 & 22 & 49 & 25994\\
9 & 16 & 21 & 56 & 25 & 48 & 1734 & 29 & 64 & 18955 & - & - & -\\
10 & 20 & 26 & 110 & 31 & 61 & 4847 & - & - & - & - & - & -\\
11 & 26 & 31 & 211 & 38 & 75 & 13532 &- & - & - & - & - & -\\
12 & 30 & 37 & 404 & 44 & 91 & 38456 & - & - & - & - & - & -\\
13 & 36 & 43 & 771 & - & - & -& - & - & - & - & - & -\\
14 & 42 & 50 & 1474 & - & - & -& - & - & - & - & - & -\\
15 & 50 & 57 & 2809 & - & - & -& - & - & - & - & - & -\\\hline
\end{tabular}
\caption{The discrepancy of the three sequences are reported for
various values of $q$ and $n$ using a Matlab code. Dashed entries
were left as they take too long.}
\end{table}

%
%

%
%


\section{Conclusion}
Relying on well-established techniques of preference functions, we
introduced two nonbinary sequences with properties that parallel the
well-known binary prefer-opposite and prefer-same sequences. The
latter two sequences were recently generated efficiently via a
simple successor rule, see Sala et al~\cite{Sala2020}. Generating
the new sequences efficiently is planned as future research. Also of
interest is an analytical study of the discrepancy of these
sequences.

\end{document}